\documentclass[11pt, reqno]{amsart}
\usepackage[french,english]{babel}
\usepackage{amssymb}
\usepackage{hyperref}
\usepackage{tikz}
\usepackage{tikz-cd}
\usepackage{array}
\usepackage[capitalise]{cleveref}
\usepackage{epsfig,graphicx}
\usepackage{float}
\usepackage{capt-of}
\tikzset{node distance=2cm, auto}
\usepackage[vcentermath]{youngtab}
\usepackage{ytableau}
\usepackage{hyperref}
\usepackage{etoolbox}
\usetikzlibrary{graphs}
\usepackage{listings}
\usetikzlibrary{backgrounds}
\usetikzlibrary{arrows,decorations.markings}
\tikzstyle{vertex}=[circle, draw, inner sep=0pt, minimum size=6pt]

\usetikzlibrary{matrix,arrows,decorations.pathmorphing}
\usepackage{mathrsfs}
\usepackage{caption}
\numberwithin{equation}{section}
\newtheorem{theorem}{Theorem}[section]
\newtheorem{corollary}[theorem]{Corollary}

\newtheorem{lemma}[theorem]{Lemma}

\newtheorem{proposition}[theorem]{Proposition}

\theoremstyle{remark}
\newtheorem{remark}[theorem]{Remark}

\title[Torus quotient of $G_{n,2n}$]{Torus quotient of the Grassmannian $G_{n,2n}$}

\author[A. Nayek]{Arpita Nayek}
\address{Arpita Nayek\\ Department of Mathematics\\ IIT Bombay, Powai, Mumbai 400076, India} 

\email{arpitan@math.iitb.ac.in}

\author[P. Saha]{Pinakinath Saha}
\address{Pinakinath Saha\\ Department of Mathematics\\ IIT Bombay, Powai, Mumbai 400076, India}

\email{psaha@math.iitb.ac.in}

\begin{document}
	
\begin{abstract}
Let $G_{n,2n}$ be the Grassmannian parameterizing the $n$-dimensional subspaces of $\mathbb{C}^{2n}.$ The Picard group of $G_{n,2n}$ is generated by a unique ample line bundle $\mathcal{O}(1).$ Let $T$ be a maximal torus of $SL(2n,\mathbb{C})$ which acts on $G_{n,2n}$ and $\mathcal{O}(1).$  By \cite[Theorem 3.10, p.764]{Kum}, $2$ is the minimal integer $k$ such that $\mathcal{O}(k)$ descends to the GIT quotient. In this article, we prove that the GIT quotient of $G_{n,2n}$ ($n\ge 3$) by $T$ with respect to $\mathcal{O}(2)=\mathcal{O}(1)^{\otimes 2}$ is not projectively normal when polarized with the descent of $\mathcal{O}(2).$ 
\end{abstract}
\maketitle

\selectlanguage{french} 
\begin{abstract}
Soit $G_{n,2n}$ la Grassmannienne des sous-espaces de dimension $n$ de $\mathbb{C}^{2n}.$ Le groupe de Picard de $G_{n,2n}$ est engendré par un unique fibré en droites ample $\mathcal{O}(1).$ Fixons un tore maximal $T$ du groupe $SL(2n,\mathbb{C})$ qui agit sur $G_{n,2n}$ et $\mathcal {O}(1).$ D'après \cite[Theorem 3.10, p.764]{Kum}, $2$ est l'entier minimal $k$ tel que $\mathcal{O}(k)$ descende au quotient GIT. Dans cet article, nous prouvons que le quotient GIT de $G_{n,2n}$ ($n\ge 3$) par $T$ par rapport à $\mathcal{O}(2)=\mathcal{O}( 1)^{\otimes 2}$ n'est pas projectivement normal lorsqu'il est polarisé avec la descente de $\mathcal{O}(2).$
\end{abstract}
\maketitle

\selectlanguage{english}
\section{Introduction}\label{section1}
A polarized variety $(X, \mathcal L),$ where $\mathcal L$ is a very ample line bundle is said to be projectively normal if its homogeneous coordinate ring $\oplus_{m \in \mathbb Z_{\geq 0}}H^0(X, \mathcal L^{\otimes m})$ is integrally closed and it is generated as a $\mathbb{C}$-algebra by $H^0(X, \mathcal L)$ (see \cite[Chapter II, Exercise 5.14]{R}). For example, the projective line $(\mathbb{P}^1,\mathcal{O}(1))$ is projectively normal. However, if we consider the rational twisted quartic curve in $\mathbb{P}^{3},$ i.e., image $X=\{[a^4:a^3b:ab^3:b^4] \in \mathbb{P}^3 : [a:b] \in \mathbb{P}^1\}$ of the embedding $i:\mathbb{P}^{1}\hookrightarrow \mathbb{P}^{3}$ given by $[a: b]\mapsto [a^{4}:a^{3}b: ab^{3}: b^{4}],$  then $(X,\mathcal{O}_{X}(1))=(\mathbb{P}^1, \mathcal{O}(3))$ is normal but not projectively normal as the affine cone of $X$ inside $\mathbb{C}^{4}$ is not normal (see \cite[Chapter I, Exercise 3.18]{R}).
	
In \cite{K1}, Kannan made an attempt to study projective normality of the GIT quotient of $G_{2,n}$ by a maximal torus $T$ of $SL(n,\mathbb{ C})$ with respect to the descent of $\mathcal{O}(n)$ ($n$ is odd). There it was proved that the homogeneous coordinate ring of the GIT quotient of $G_{2,n}$ by $T$ with respect to the descent of $\mathcal{O}(n)$ is a finite module over the subring generated by the degree one elements. In \cite{howard}, Howard et al. showed that the GIT quotient of $G_{2,n}$ by $T$ with respect to the descent of $\mathcal{O}(\frac{n}{2})$ (respectively, $\mathcal{O}(n)$) is projectively normal if $n$ is even (respectively, if $n$ is odd). In \cite{NP}, Nayek et al. used graph theoretic techniques to give a short proof of the projective normality of the GIT quotient of $G_{2,n}$ by $T$ with respect to the descent of $\mathcal{O}(n)$ for any $n.$
	
To the best of our knowledge it is not known whether there is a suitable ample line bundle $\mathcal{L}$ on $G_{r,n}$ ($r\ge 3$) such that the GIT quotient of $G_{r,n}$ by $T$ with respect to the descent of the line bundle $\mathcal{L}$ is projectively normal (respectively, not projectively normal) with respect to the descent of $\mathcal{L}.$ 
	
In this article, we prove the following:
	
\begin{theorem}\label{thm1}
The GIT quotient of $G_{n,2n}$ $(n \geq 3)$ by a maximal torus $T$ of $SL(2n,\mathbb{ C})$ with respect to the descent of $\mathcal{O}(2)$ is not projectively normal $($for more precise see \cref{corollary3.5}$).$ 
\end{theorem}

The layout of the paper is as follows. In \cref{section2}, we recall some preliminaries on algebraic groups, Standard Monomial Theory and Geometric Invariant Theory.  In \cref{section3}, we prove \cref{thm1} (see \cref{corollary3.5}).

	\section{Notation and Preliminaries}\label{section2}
    We refer to \cite{Hum1}, \cite{Hum2},\cite{LB}, \cite{Mumford}, \cite{New} and \cite{Ses} for preliminaries in algebraic groups, Lie algebras, Standard Monomial Theory and Geometric Invariant Theory. 
	
	Let $V=\mathbb{C}^{2n}$ and $(e_{1},e_{2},\ldots,e_{2n})$ be the standard basis of $V.$ For a fixed integer $r$ with $1\le r\le 2n-1,$ let $G_{r,2n}$ be the Grassmannian parameterizing the $r$-dimensional subspaces of $\mathbb{C}^{2n}.$ Then there is a natural projective variety structure on $G_{r,2n}$ given by the Pl\"ucker embedding $\pi: G_{r,2n}\hookrightarrow \mathbb{P}(\wedge^{r}V)$ sending $r$-dimensional subspace to its $r$-th exterior power. The natural left action of $SL(2n,\mathbb{ C})$ on $V$ induces an action of $SL(2n,\mathbb{ C})$ on $\wedge^{r}V$ and thus on $\mathbb{P}(\wedge^{r}V),$ moreover, $\pi$ is $SL(2n, \mathbb{C})$-equivariant. Let $T$ be the maximal torus of $SL(2n,\mathbb{C})$ consisting of diagonal matrices. Let $\mathcal{O}(1)$ denote the hyperplane line bundle on $G_{r,2n}$ given by the Pl\"{u}cker embedding $\pi.$ Note that $\mathcal{O}(1)$ is $SL(2n, \mathbb{C})$-linearized, in particular, $T$-linearized. 
	
	 Let $I(r,2n)$ denote the indexing set $\{\underline{i}=(i_{1},i_{2},\ldots, i_{r})| i_{j}\in \mathbb{Z} \text{~and~} 1\le i_{1}<i_{2}<\cdots <i_{r}\le 2n\}.$ Let $e_{\underline{i}}=e_{i_1} \wedge e_{i_2} \wedge \cdots \wedge e_{i_r}$ for $\underline{i}=(i_1, i_2, \ldots, i_r) \in I(r,2n).$ Then $\{e_{\underline{i}}: \underline{i} \in I(r,2n)\}$ forms a basis of $\wedge^rV.$ Let $\{p_{\underline{i}}: \underline{i}\in I(r,2n)\}$ be the basis of the dual space $(\wedge^{r}V)^{*},$ which is dual to $\{e_{\underline{i}}: \underline{i} \in I(r,2n)\},$ i.e., $p_{\underline{j}}(e_{\underline{i}})=\delta_{ij}.$ Note that $p_{\underline{i}}$'s are the $\underline{i}^{\text{th}}$ Pl\"{u}cker coordinates of $G_{r,2n}.$ 
	 
	 In $V,$ we fix a full flag $\{0\}=V_{0}\subset V_{1}\subset \cdots\subset V_{2n}=V.$ For $w=(w_1,w_2,\ldots,w_r)$ in $I(r,2n),$ the Schubert variety in $G_{r,2n}$  associated to $w$ is denoted by $X(w)$ and is defined by $$X(w)=\bigg\{W\in G_{r, 2n}| \begin{matrix}
	 	\dim W\cap V_{j}\ge i, \text{~if~} w_i\le j< w_{i+1},\\ \text{~where~} 1\le j\le 2n, 0\le i\le r \text{~and~} w_{0}:=0 ,w_{r+1}:=2n
	 \end{matrix}\bigg\}.$$

	 The definition of a Schubert variety $X(w)$ depends on the choice of a full flag. However, given any two full flags $\{0\}=V_{0}\subset V_{1}\subset \cdots\subset V_{2n}=V$ and $\{0\}=V'_{0}\subset V'_{1}\subset \cdots\subset V'_{2n}=V$ in $V,$ there exist an automorphism of $V$ which takes $V_i$ to $V'_{i},$ which shows that $X(w)$ is well defined up to an automorphism of $V.$ We note that $X(w)$ is a closed subvariety of $G_{r,2n}$ of dimension $\displaystyle\sum_{i=1}^{r} w_{i}-\frac{r(r+1)}{2}.$  
	 
	 	There is a natural partial order on $I(r,2n),$ given as follows: for $v=(v_{1},v_{2},\ldots, v_{r}),$ $w=(w_{1},w_{2},\ldots, w_{r}),$ $v\le w$ if and only if $v_{i} \leq w_{i}$ for all $1\le i\le r.$ For $v, w \in I(r,2n),$ $X(v) \subseteq X(w)$ if and only if $v \leq w.$ Further, $p_{v}|_{X(w)} \neq 0$ if and only if $v \leq w.$ 
	 	
	 	For $w\in I(r,2n),$ we also denote the restriction of the line bundle  $\mathcal{O}(1)$ on $G_{r,2n}$ to $X(w)$ by $\mathcal{O}(1).$ The monomial $p_{\tau_1}p_{\tau_2}\ldots p_{\tau_m} \in H^0(X(w),\mathcal{O}(m)),$ where $\tau_1, \tau_2, \ldots, \tau_m \in I(r,2n)$ is said to be standard monomial of degree $m$ if $\tau_1\leq \tau_2\leq \cdots \leq \tau_m \leq w.$ The standard monomials of degree $m$ on $X(w)$ form a basis of $H^0(X(w),\mathcal{O}(m)).$
	The Grassmannian $G_{r,2n} \subseteq \mathbb{P}(\wedge^r V)$ is precisely the zero set of the following well known Pl\"{u}cker relations:
	\begin{align}
		\sum_{h=1}^{r+1}(-1)^h p_{i_1,i_2,\ldots,i_{r-1}j_h}p_{j_1,\ldots,\hat{j_h},\ldots,j_{r+1}}
		\label{2.1},
	\end{align}
	where $\{i_1,\ldots, i_{r-1}\},$ $\{j_1,\ldots,j_{r+1}\}$ are two subsets of $\{1,2,\ldots,2n\}$ and $\hat{j_{h}}$ means dropping the index $j_{h}.$
	
	  A point $p \in X(w)$ is said to be semi-stable with respect to the $T$-linearized line bundle $\mathcal{O}(1)$ if there is a $T$-invariant section $s \in H^0(X(w),\mathcal{O}(m))$ for some positive integer $m$ such that $s(p)\neq 0.$ We denote the set of all semi-stable points of $X(w)$ with respect to $\mathcal{O}(1)$ by $X(w)^{ss}_T(\mathcal{O}(1)).$ A point $p$ in $X(w)^{ss}_{T}(\mathcal{O}(1))$ is said to be stable if the $T$-orbit of $p$ is closed in $X(w)^{ss}_{T}(\mathcal{O}(1))$ and the stabilizer of $p$ in $T$ is finite. We denote the set of all stable points of $X(w)$ with respect to $\mathcal{O}(1)$ by $X(w)^{s}_T(\mathcal{O}(1)).$ 	
	
	Let $B~(\supset T)$ be the Borel subgroup of $SL(2n,\mathbb{ C})$ consisting of upper triangular matrices. For $1\le i\le 2n,$ define $\varepsilon_{i}:T\to \mathbb{C}^{\times}$ by  $\varepsilon_{i}(\text{diag}(t_{1},\ldots ,t_{2n}))=t_{i}.$ Then  $S:=\{\alpha_i:=\varepsilon_{i}-\varepsilon_{i+1}| \text{~for all ~} 1\le i\le 2n-1\}$ forms the set of simple roots of $SL(2n, \mathbb{C})$ with respect to $T$ and $B.$ Let $\{\varpi_i| i= 1,2, \ldots, 2n-1\}$ be the set of fundamental dominant weights corresponding to $S.$
	 
	 	For $\lambda=m\varpi_r$ $(m \geq 1)$, we associate a Young diagram (denoted by $\Gamma$) with $\lambda_i$ number of boxes in the $i$-th column, where $\lambda_i:=m$ for $1 \leq i \leq r$. It is also called Young diagram of shape $\lambda.$	A Young diagram $\Gamma$ associated to $\lambda$ is said to be a Young tableau if the diagram is filled with integers $1, 2, \ldots, 2n$. We also denote this Young tableau by $\Gamma$. A Young tableau is said to be standard if the entries along any column is non-decreasing from top to bottom and along any row is strictly increasing from left to right. Given a Young tableau $\Gamma$, let $\tau=\{i_1,i_2,\ldots,i_r\}$ be a typical row in $\Gamma$, where $1 \leq i_1 < \cdots < i_r \leq 2n$. To the row $\tau$, we associate the Pl\"{u}cker coordinate $p_{i_1,i_2, \ldots,i_r}$. We set $p_{\Gamma}=\prod_{\tau}p_{\tau}$, where the product is taken over all the rows of $\Gamma$. Note that for $w\in I(r,2n),$ $p_{\Gamma}$ is a standard monomial on $X(w)$ if $\Gamma$ is standard and the bottom row of $\Gamma$ is less than or equal to $w.$ Further,  $p_{\Gamma}$ is also called standard monomial on $X(w)$ of shape $\lambda.$ We use the notation $p_{\Gamma}$ and $\Gamma$ interchangeably.
	
	 Now we recall the definition of weight of a standard Young tableau $\Gamma$ (see \cite[Section 2, p.336]{L}). For a positive integer $1\leq i\leq 2n$, we denote by $c_{\Gamma}(i)$,
	the number of boxes of $\Gamma$ containing the integer $i$. 
	The weight of $\Gamma$ is defined as 
	$wt(\Gamma):=c_{\Gamma}(1)\varepsilon_1+ \cdots + c_{\Gamma}(2n)\varepsilon_{2n}.$ \noindent 
	We conclude this section by recalling the following key lemma about $T$-invariant monomials in $H^0(G_{r,2n}, \mathcal{O}(m)).$
	\begin{lemma}$($See \cite[Lemma 3.1, p.4]{NP}$)$\label{lemma2.1}
		A monomial $p_{\Gamma}\in H^0(G_{r,2n}, \mathcal{O}(m))$ is $T$-invariant if and only if  $c_{\Gamma}(i)=c_{\Gamma}(j)$ for all $1 \leq i,j \leq 2n.$
	\end{lemma}

	\section{Main Theorem} \label{section3}
	First we recall that by \cite[Theorem 3.10, p.764]{Kum}, $2$ is the minimal integer $k$ such that the line bundle $\mathcal{O}(k)$ on $G_{n,2n}$ descends to the GIT quotient $T \backslash \backslash (G_{n,2n})^{ss}_T(\mathcal{O}(2)).$ In this section, we prove that there exists a Schubert subvariety $X(v)$ of $G_{n,2n}$ admitting semi-stable points such that the GIT quotient $T\backslash\backslash(X(v))^{ss}_T(\mathcal{O}(2))$ with respect to the descent of $\mathcal{O}(2)$ is not projectively normal (see \cref{theorem3.4}). As a consequence, we conclude that any Schubert variety $X(w)$ containing $X(v),$ the GIT quotient $T\backslash\backslash(X(w))^{ss}_T(\mathcal{O}(2))$ with respect to the descent of $\mathcal{O}(2)$ is not projectively normal. In particular, $T\backslash\backslash(G_{n,2n})^{ss}_T(\mathcal{O}(2))$ is not projectively normal.
	
	Recall that $2\varpi_n=2\varepsilon_1+2\varepsilon_2+\cdots+2\varepsilon_n$ (see \cite[Table 1, p.69]{Hum1}). For $u=(u_1,u_2,\ldots,u_n) \in I(n,2n),$ we define $u(2\varpi_n)=2\varepsilon_{u_{1}}+2\varepsilon_{u_2}+\cdots+2\varepsilon_{u_n}.$   
	
	Recall that by \cite[Corollary 1.9, p.85]{KS}, there exists a unique minimal element $w_1 \in I(n,2n)$ such that $w_1(2\varpi_{n})\leq 0,$ i.e., $-w_1(2\varpi_n)$ is a non-negative linear combination of simple roots. Consider $w=(2,4,6,\ldots,2n-4,2n-2,2n).$ Then $-w(2\varpi_n)=-(2\varepsilon_2+2\varepsilon_4+\cdots+2\varepsilon_{2n})=\alpha_1+\alpha_3+\cdots+\alpha_{2n-1},$ as $\sum_{i=1}^{2n}\varepsilon_{i}=0.$ Thus, $w(2\varpi_n) \leq 0.$ On the other hand for any $v \le w$ such that $\sum_{i=1}^{n}({w}_i-v_i)=1,$ we have $v=(2,4,6,\ldots,2i-2,2i-1,2i+2,\ldots,2n-4,2n-2,2n)$ for some $1 \leq i \leq n.$ Then $-v(2\varpi_n)=-(2\varepsilon_2+\cdots+2\varepsilon_{2i-2}+2\varepsilon_{2i-1}+2\varepsilon_{2i+2}+\cdots+2\varepsilon_{2n}).$ Since $\sum_{i=1}^{2n}\varepsilon_{i}=0,$ we have $-v(2\varpi_n)=(\alpha_1+\alpha_3+\cdots+\alpha_{2i-3}+\alpha_{2i+1}+\alpha_{2i+3}+\cdots+\alpha_{2n-1})-\alpha_{2i-1}.$ Thus, $v(2\varpi_{n}) \nleq 0.$ Therefore, $w=w_1.$
	
	Now we consider the following $w_i$'s such that $w_1 \leq w_i$ for all $2 \leq i \leq 5:$
   \begin{itemize} 
		\item $w_2~~=~(2,4,6,\ldots,2n-6,2n-3,2n-2,2n)$
		\item $w_3~~=~(2,4,6,\ldots,2n-6,2n-4,2n-1,2n)$
		\item $w_4~~=~(2,4,6,\ldots,2n-6,2n-3,2n-1,2n)$
		\item  $w_5~~=~(2,4,6,\ldots,2n-6,2n-2,2n-1,2n).$
	\end{itemize}
	
Note that $\{ w_1, w_2, w_3,w_4, w_5 \}$ is precisely the set $\{ w\in I(n,2n): w_1 \leq w \leq w_5 \}.$ Further, note that $w_{2}$ and $w_{3}$ are non-comparable and $w_{2},w_{3}\le w_{4}\le w_{5}.$
Since $w_1 \leq w_{i}$ and $w_{1}(2\varpi_{n})\leq 0,$ we have $w_{i}(2\varpi_{n}) \leq 0$ for all $2 \leq i \leq 5.$ Thus, by \cite[Lemma 2.1, p.470]{KP}, $X(w_i)^{ss}_T(\mathcal{O}(2))$ is non-empty for all $1 \leq i \leq 5.$ 

Let $X= T \backslash \backslash (X(w_{5}))^{ss}_T(\mathcal{O}(2)).$
Then we have $X= {\rm Proj}(R),$ where $R=\bigoplus\limits_{k\in\mathbb{Z}_{\geq 0}} R_{k}$ and $R_{k}=H^0(X(w_{5}), \mathcal{O}(2k))^{T}.$ Note that $R_{k}$'s  are finite dimensional vector spaces. 

Let us consider the following standard monomials
	
	\tiny{\ytableausetup{boxsize=3.2em} \[X_1= \ytableausetup{centertableaux}
		\begin{ytableau}
			1 & 3 & 5  & \cdots & 2n-7 & 2n-5 & 2n-3 & 2n-1 \\
			2 & 4 & 6 & \cdots & 2n-6 & 2n-4 & 2n-2 & 2n  
		\end{ytableau} 
		~ X_2= \ytableausetup{centertableaux}
		\begin{ytableau}
			1 & 3 & 5  & \cdots & 2n-7 & 2n-5 & 2n-4 & 2n-1 \\
			2 & 4 & 6 & \cdots & 2n-6 & 2n-3 & 2n-2 & 2n  
		\end{ytableau}\] 
		\[X_3= \ytableausetup{centertableaux}
		\begin{ytableau}
			1 & 3 & 5  & \cdots & 2n-7 & 2n-5 & 2n-3 & 2n-2 \\
			2 & 4 & 6 & \cdots & 2n-6 & 2n-4 & 2n-1 & 2n  
		\end{ytableau} 
		~	X_4= \ytableausetup{centertableaux}
		\begin{ytableau}
			1 & 3 & 5  & \cdots & 2n-7 & 2n-5 & 2n-4 & 2n-2 \\
			2 & 4 & 6 & \cdots & 2n-6 & 2n-3 & 2n-1 & 2n  
		\end{ytableau}\] 
		$X_5= \ytableausetup{centertableaux}
		\begin{ytableau}
			1 & 3 & 5  & \cdots & 2n-7 & 2n-5 & 2n-4 & 2n-3 \\
			2 & 4 & 6 & \cdots & 2n-6 & 2n-2 & 2n-1 & 2n  
		\end{ytableau}$ }
	
	\tiny{\ytableausetup{boxsize=3.2em} \[Y_1= \ytableausetup{centertableaux}
		\begin{ytableau}
			1& 3 & 5  & \cdots & 2n-7 & 2n-5 & 2n-4 & 2n-3\\
			1& 3 & 5  & \cdots & 2n-7 & 2n-5 & 2n-2 & 2n-1\\
			2 & 4 & 6 & \cdots & 2n-6 & 2n-4 & 2n-2 & 2n\\
			2 & 4 & 6 & \cdots & 2n-6 & 2n-3 & 2n-1 & 2n 
		\end{ytableau}
		~Y_2= \ytableausetup{centertableaux}
		\begin{ytableau}
			1& 3 & 5  & \cdots & 2n-7 & 2n-5 & 2n-4 & 2n-2\\
			1& 3 & 5  & \cdots & 2n-7 & 2n-5 & 2n-3 & 2n-1\\
			2 & 4 & 6 & \cdots & 2n-6 & 2n-4 & 2n-3 & 2n\\
			2 & 4 & 6 & \cdots & 2n-6 & 2n-2 & 2n-1 & 2n 
		\end{ytableau}.\]}
	\normalsize
	
	\begin{remark}\label{rem3.6}
		Let $\mathcal{M}$ denote the descent of the line bundle $\mathcal{O}(2)$ to $X.$ Then by using {\it Quantization commutes with reduction}  we have $H^0(X(w_{5}), \mathcal{O}(2k))^{T}=H^{0}(X,\mathcal{M}^{\otimes k})$ for $k\in \mathbb{Z}_{\ge 0}$ (see \cite[Theorem 3.2.a., p.11]{Tel} or \cite[Theorem 4.1(ii), p.526]{Ses1}).
	\end{remark}
			\begin{remark}\label{remark3.2}
		\begin{itemize}
			\item[(i)] Note that the set of $T$-invariant standard monomials of shape $2\varpi_{n}$ on $X(w_5)$ is $\{X_{i}$: $1\leq i \leq 5$\}. Thus by \cite[Theorem 12.4.8, p.207]{LB}, the set $\{X_i: 1\leq i \leq 5\}$ forms standard monomial basis of $R_1.$ 
			
	\item[(ii)]	The set of $T$-invariant standard monomials of shape $4\varpi_{n}$ on $X(w_5)$ is $\{Y_{1}, Y_{2}\} \cup \{X_{i}X_{j} :1\le i\le j\le 5\}\setminus \{X_2X_3\}.$ Therefore, by \cite[Theorem 12.4.8, p.207]{LB}, the set $\{Y_{1}, Y_{2}\}\cup \{X_{i}X_{j}: 1\le i\le j\le 5 \}\setminus\{X_{2}X_{3}\}$ forms standard monomial basis of $R_{2}.$
		\end{itemize}
	\end{remark}

	\begin{theorem}\label{theorem3.4}
	The GIT quotient $X$ with respect to the descent of $\mathcal{O}(2)$ is not projectively normal.	
\end{theorem}
\begin{proof}
	 Consider the natural map $f: R_{1}\otimes R_{1}\to R_{2}$ of vector spaces given by $X_{i}\otimes X_{j}\mapsto X_{i}X_{j}$ for $1\le i,j\le 5.$ Then $f$ factors through second symmetric power $S^{2}R_{1}$ of the vector space $R_{1}.$ For simplicity we also denote the factor map $S^{2}R_{1}\to R_{2}$ by $f.$ 
	By \cref{remark3.2}, we have $\dim(R_{1})=5$ and $\dim(R_{2})=16.$ So, the map $f: S^{2}R_{1}\to R_{2}$ cannot be surjective, since $\dim(S^{2}R_{1})=15<\dim(R_{2}).$ 
\end{proof}
\begin{corollary}\label{corollary3.5}
	The GIT quotient $T \backslash \backslash (X(w))^{ss}_T(\mathcal{O}(2))$ with respect to the descent of $\mathcal{O}(2)$ is not projectively normal for $w \in I(n,2n)$ such that $w_{5} \leq w$. In particular, the GIT quotient $T \backslash \backslash (G_{n,2n})^{ss}_T(\mathcal{O}(2))$ with respect to the descent of $\mathcal{O}(2)$ is not projectively normal.
\end{corollary}
\begin{proof}
	By \cite[Theorem 3.1.1(b), p.85]{BK}, the restriction map $$\phi: H^{0}(X(w), \mathcal{O}(2k)) \to H^{0}(X(w_{5}), \mathcal{O}(2k))$$ is surjective. Further, since $T$ is linearly reductive, the restriction map $\phi: H^{0}(X(w), \mathcal{O}(2k))^T \to H^{0}(X(w_{5}), \mathcal{O}(2k))^T$ is surjective for all $k \geq 1$. So, by \cref{theorem3.4}, $T \backslash \backslash (X(w))^{ss}_T(\mathcal{O}(2))$ with respect to the descent of $\mathcal{O}(2)$ is not projectively normal.
\end{proof}
	\begin{lemma}\label{lemma3.1}
			The homogeneous coordinate ring of $X$ is generated by elements of degree at most two.
	\end{lemma}
	\begin{proof}
		Let $f \in R_k$ be a standard monomial. We claim that $f = f_1f_2,$ where $f_1$ is in $R_1$ or $R_2$.
		The Young diagram associated to $f$ has the shape $(\lambda_1,  \lambda_2, \lambda_3, \ldots, \lambda_{n}) = (\underbrace{2k, 2k, \ldots, 2k}_{n})$. So the Young tableau $\Gamma$ associated to $f$ has $2k$ rows and $n$ columns with strictly increasing rows and non-decreasing columns. Since $f$ is $T$-invariant, by \cref{lemma2.1}, we have $c_\Gamma(t)=k$ for all $1 \leq t \leq 2n.$
		Let $r_i$ be the $i$-th row of the tableau. Let $E_{i,j}$ be the $(i,j)$-th entry of the tableau $\Gamma$ and $N_{t,j}$ is the number of boxes in the $j$-th column of $\Gamma$ containing the integer $t.$ 
		
		Recall that $w_5$ is $(2,4,6,\ldots,2n-6,2n-2,2n-1,2n).$ Since $r_{2k} \leq w_5,$ we have $E_{2k,j} \leq 2j$ for all $1 \leq j \leq n-3.$ Note that for $1 \leq j \leq n-3,$ we have $E_{i,j}=2j$ for all $k+1 \leq i \leq 2k.$ Thus, for $1 \leq j \leq n-2,$ we have $E_{i,j}=2j-1$ for all $1 \leq i \leq k.$ Thus, the following rows are the possibilities for $r_{2k}:$
		\begin{itemize}
			\item $(2,4,6,\ldots,2n-6,2n-4,2n-2,2n)$
			\item $(2,4,6,\ldots,2n-6,2n-3,2n-2,2n)$
			\item $(2,4,6,\ldots,2n-6,2n-4,2n-1,2n)$
			\item $(2,4,6,\ldots,2n-6,2n-3,2n-1,2n)$
			\item  $(2,4,6,\ldots,2n-6,2n-2,2n-1,2n).$
		\end{itemize} 
		\textbf{Case I:} Assume that $r_{2k}=(2,4,6,\ldots,2n-6,2n-4,2n-2,2n).$ Then for $1 \leq j \leq n,$ we have $E_{i,j}=2j-1$ (resp. $E_{i,j}=2j$) for all $1 \leq i \leq k$ (resp. for all $k+1 \leq i \leq 2k$). Therefore, $r_1=(1,3,5,\ldots,2n-7,2n-5,2n-3,2n-1).$ Hence, $r_1, r_{2k}$ together give a factor $X_{1}$ of $f.$
		
		\textbf{Case II:} Assume that $r_{2k}=(2,4,6,\ldots,2n-6,2n-3,2n-2,2n).$ Since $N_{2n-3,n-2}\geq 1,$ we have $N_{2n-4,n-2} \leq k-1.$ Thus, $E_{1,n-1}=2n-4.$ Since $E_{2k,n-1}=2n-2,$ we have $E_{i,n} = 2n-1$ for all $1 \leq i \leq k.$ Therefore, $r_1=(1,3,5,\ldots,2n-7,2n-5,2n-4,2n-1).$ Hence, $r_1, r_{2k}$ together give a factor $X_{2}$ of $f.$  
		
		\textbf{Case III:} Assume that $r_{2k}=(2,4,6,\ldots,2n-6,2n-4,2n-1,2n).$ Then $E_{i,n-2}=2n-4$ for all $k+1 \leq i \leq 2k$ and $E_{i,n-1}=2n-3$ for all $1 \leq i \leq k.$ Since $N_{2n-1,n-1}\geq 1,$ we have $N_{2n-2,n-1} \leq k-1.$ Thus, $E_{1,n}=2n-2.$ Therefore, $r_1=(1,3,5,\ldots,2n-7,2n-5,2n-3,2n-2).$ Hence, $r_1, r_{2k}$ together give a factor $X_{3}$ of $f.$ 
		
		\textbf{Case IV:} Assume that $r_{2k}=(2,4,6,\ldots,2n-6,2n-3,2n-1,2n).$ Since $N_{2n-3,n-2}\geq 1,$ we have $N_{2n-4,n-2} \leq k-1.$ Thus, $E_{1,n-1}=2n-4.$ Since $E_{2k,n-1}=2n-1,$ we have $E_{1,n} \leq 2n-2.$ Thus, $r_1$ is either $(1,3,5,\ldots,2n-7,2n-5,2n-4,2n-3)$ or $(1,3,5,\ldots,2n-7,2n-5,2n-4,2n-2).$  If $r_1=(1,3,5,\ldots,2n-7,2n-5,2n-4,2n-2),$ 
		then $r_1, r_{2k}$ together give a factor $X_{4}$ of $f.$ If $r_1=(1,3,5,\ldots,2n-7,2n-5,2n-4,2n-3),$ then $E_{k,n-1}=2n-2.$ Otherwise, if $E_{k,n-1}=2n-4$ or $2n-3,$ then $\sum_{t=2n-4}^{2n-3}(N_{t, n-2}+N_{t, n-1}+N_{t,n})\ge 2k+1,$ which is a contradiction. Therefore, $E_{k,n}=2n-1.$ Thus, $E_{k+1,n-1}=2n-2.$ Therefore, $r_{k}=(1,3,5,\ldots,2n-7,2n-5,2n-2,2n-1)$ and $r_{k+1}=(2,4,6,\ldots,2n-6,2n-4,2n-2,2n).$ Hence, $r_1, r_{k}, r_{k+1}, r_{2k}$ together give a factor $Y_1$ of $f.$
		
		\textbf{Case V:} Assume that $r_{2k}=(2,4,6,\ldots,2n-6,2n-2,2n-1,2n).$ Since $E_{2k,n-2}=2n-2,$ we have $E_{1,n-1}=2n-4.$  Thus, $E_{1,n} \geq 2n-3.$ If $E_{1,n} = 2n-1,$ then $N_{2n-1,n-1}+N_{2n-1,n} \geq k+1,$ which is a contradiction. Thus, $E_{1,n}$ is either $2n-3$ or $2n-2.$ Thus, $r_1$ is either $(1,3,5,\ldots,2n-7,2n-5,2n-4,2n-3)$ or $(1,3,5,\ldots,2n-7,2n-5,2n-4,2n-2).$  If $r_1=(1,3,5,\ldots,2n-7,2n-5,2n-4,2n-3),$ then $r_1, r_{2k}$ together give a factor $X_{5}$ of $f.$ If $r_1=(1,3,5,\ldots,2n-7,2n-5,2n-4,2n-2),$ then $2n-4 \leq E_{k,n-1} \leq 2n-1.$ If $E_{k,n-1}=2n-1,$ then $N_{2n-1,n-1} \geq k+1,$ which is a contradiction. If $E_{k,n-1}=2n-2,$ then $\sum_{t=2n-2}^{2n-1}(N_{t,n-2}+N_{t,n-1}+N_{t,n}) \geq 2k+2,$ which is a contradiction. If $E_{k,n-1}=2n-4,$ then $E_{i,n-2}=2n-3$ for all $k+1 \leq i \leq 2k-1.$ Thus, $c_{\Gamma}(2n-3) \leq k-1,$ which is a contradiction. Therefore, $E_{k,n-1}=2n-3.$ Thus, $E_{k,n}=2n-1.$ Then $2n-3 \leq E_{k+1,n-1} \leq 2n-1.$ If $E_{k+1,n-1}=2n-1,$ then $N_{2n-1,n-1}+N_{2n-1,n} \geq k+1,$ which is a contradiction. If $E_{k+1,n-1}=2n-2,$ then $\sum_{t=2n-2}^{2n-1}(N_{t,n-2}+N_{t,n-1}+N_{t,n}) \geq 2k+1,$ which is a contradiction. Thus, $E_{k+1,n-1}=2n-3.$ Therefore, $r_{k}=(1,3,5,\ldots,2n-7,2n-5,2n-3,2n-1)$ and $r_{k+1}=(2,4,6,\ldots,2n-6,2n-4,2n-3,2n).$ Therefore, $r_1, r_{k}, r_{k+1}, r_{2k}$ together give a factor $Y_2$ of $f.$
	\end{proof}
	\begin{lemma}\label{lemma3.3}
		$X_i$'s $(1 \leq i \leq 5),$ and $Y_j$'s $(1 \leq j \leq 2)$ satisfy the following relation in $R_2:$
		$X_2X_3=X_1X_4-Y_2-Y_1+X_5(X_1-X_2-X_3+X_4-X_5).$
	\end{lemma}
	\begin{proof} 
		Note that	\tiny{\ytableausetup{boxsize=3.2em} \ytableausetup{centertableaux}
			\[X_2X_3=\begin{ytableau}
				1 & 3 & 5  & \cdots & 2n-7 & 2n-5 & 2n-4 & 2n-1 \\
				1 & 3 & 5  & \cdots & 2n-7 & 2n-5 & 2n-3 & 2n-2 \\
				2 & 4 & 6 & \cdots & 2n-6 & 2n-3 & 2n-2 & 2n \\
				2 & 4 & 6 & \cdots & 2n-6 & 2n-4 & 2n-1 & 2n
			\end{ytableau}.\]}
		\normalsize{By} using \eqref{2.1}, we have the following straightening laws in $X(w_5)$
		
		\tiny{\ytableausetup{boxsize=3.2em} \ytableausetup{centertableaux}
			\[\begin{ytableau}
				1 & 3 & 5  & \cdots & 2n-7 & 2n-5 & 2n-4 & 2n-1 \\
				1 & 3 & 5  & \cdots & 2n-7 & 2n-5 & 2n-3 & 2n-2
			\end{ytableau}\]
			\[=\hspace{.15cm}\begin{ytableau}
				1 & 3 & 5  & \cdots & 2n-7 & 2n-5 & 2n-4 & 2n-2 \\
				1 & 3 & 5  & \cdots & 2n-7 & 2n-5 & 2n-3 & 2n-1
			\end{ytableau}
			~-\hspace{.15cm}\begin{ytableau}
				1 & 3 & 5  & \cdots & 2n-7 & 2n-5 & 2n-4 & 2n-3 \\
				1 & 3 & 5  & \cdots & 2n-7 & 2n-5 & 2n-2 & 2n-1
			\end{ytableau}\]}
		
		\normalsize{and} 
		
		\tiny{\ytableausetup{boxsize=3.2em} \ytableausetup{centertableaux}
			\[\begin{ytableau}
				2 & 4 & 6 & \cdots & 2n-6 & 2n-3 & 2n-2 & 2n \\
				2 & 4 & 6 & \cdots & 2n-6 & 2n-4 & 2n-1 & 2n 
			\end{ytableau}\]
			\[ = \hspace{.15cm}\begin{ytableau}
				2 & 4 & 6 & \cdots & 2n-6 & 2n-4 & 2n-2 & 2n \\
				2 & 4 & 6 & \cdots & 2n-6 & 2n-3 & 2n-1 & 2n
			\end{ytableau}
			~-\hspace{.15cm} \hspace{.15cm}\begin{ytableau}
				2 & 4 & 6 & \cdots & 2n-6 & 2n-4 & 2n-3 & 2n \\
				2 & 4 & 6 & \cdots & 2n-6 & 2n-2 & 2n-1 & 2n
			\end{ytableau}.\]
		}
			
		\normalsize{Therefore,} by using the above straightening laws we have
		
		$X_2X_3=$
		
		\tiny{\ytableausetup{boxsize=3.2em} \ytableausetup{centertableaux}
			\hspace{.4cm}\begin{ytableau}
				1 & 3 & 5  & \cdots & 2n-7 & 2n-5 & 2n-4 & 2n-2 \\
				1 & 3 & 5  & \cdots & 2n-7 & 2n-5 & 2n-3 & 2n-1 \\
				2 & 4 & 6 & \cdots & 2n-6 & 2n-4 & 2n-2 & 2n \\
				2 & 4 & 6 & \cdots & 2n-6 & 2n-3 & 2n-1 & 2n
			\end{ytableau}
			~$-$~\hspace{.15cm}\hspace{.15cm}\begin{ytableau}
				1 & 3 & 5  & \cdots & 2n-7 & 2n-5 & 2n-4 & 2n-2 \\
				1 & 3 & 5  & \cdots & 2n-7 & 2n-5 & 2n-3 & 2n-1 \\
				2 & 4 & 6 & \cdots & 2n-6 & 2n-4 & 2n-3 & 2n \\
				2 & 4 & 6 & \cdots & 2n-6 & 2n-2 & 2n-1 & 2n 
			\end{ytableau}\\
			
			\vspace{.1cm}
			~$-$~\begin{ytableau}
				1 & 3 & 5  & \cdots & 2n-7 & 2n-5 & 2n-4 & 2n-3 \\
				1 & 3 & 5  & \cdots & 2n-7 & 2n-5 & 2n-2 & 2n-1 \\
				2 & 4 & 6 & \cdots & 2n-6 & 2n-4 & 2n-2 & 2n \\
				2 & 4 & 6 & \cdots & 2n-6 & 2n-3 & 2n-1 & 2n
			\end{ytableau}
			~$+$ ~\hspace{.15cm}\hspace{.15cm}\begin{ytableau}
				1 & 3 & 5  & \cdots & 2n-7 & 2n-5 & 2n-4 & 2n-3 \\
				1 & 3 & 5  & \cdots & 2n-7 & 2n-5 & 2n-2 & 2n-1 \\
				2 & 4 & 6 & \cdots & 2n-6 & 2n-4 & 2n-3 & 2n \\
				2 & 4 & 6 & \cdots & 2n-6 & 2n-2 & 2n-1 & 2n
		\end{ytableau}}
		
		\normalsize{$=X_1X_4-Y_2-Y_1+X_5Z$}, where $Z=$\tiny{\ytableausetup{boxsize=3.2em} \ytableausetup{centertableaux}\begin{ytableau}
				1 & 3 & 5  & \cdots & 2n-7 & 2n-5 & 2n-2 & 2n-1 \\
				2 & 4 & 6 & \cdots & 2n-6 & 2n-4 & 2n-3 & 2n 
			\end{ytableau}.}
		
	\normalsize{By} using $(*)$ (see Appendix) we have $ Z = X_1-X_2-X_3+X_4-X_5.$
		Therefore, $X_2X_3 = X_1X_4-Y_2-Y_1+X_5(X_1-X_2-X_3+X_4-X_5).$
	\end{proof}

	\begin{proposition} We have
		\begin{itemize}
			\item[(i)] The GIT quotient $T\backslash\backslash(X(w_1))^{ss}_T(\mathcal{O}(2))$ with respect to the descent of $\mathcal{O}(2)$ is projectively normal and isomorphic to point.
			\item[(ii)] The GIT quotient $T\backslash\backslash(X(w_2))^{ss}_T(\mathcal{O}(2))$ with respect to the descent of $\mathcal{O}(2)$ is projectively normal and isomorphic to $\mathbb{P}^{1}$ polarized with $\mathcal{O}(1).$ 
			\item[(iii)]The GIT quotient $T\backslash\backslash(X(w_3))^{ss}_T(\mathcal{O}(2))$  with respect to the descent of $\mathcal{O}(2)$ is projectively normal and isomorphic to $\mathbb{P}^{1}$ polarized with $\mathcal{O}(1).$
			\item[(iv)] The GIT quotient $T\backslash\backslash(X(w_4))^{ss}_T(\mathcal{O}(2))$ with respect to the descent of $\mathcal{O}(2)$ is projectively normal and isomorphic to $\mathbb{P}^{3}$ polarized with $\mathcal{O}(1).$
		\end{itemize}	
	\end{proposition}
	\begin{proof} Note that $Y_{2}=0$ on $X(w_{i})$ for all $1 \leq i \leq 4.$ 
		
		Proof of (iv): By \cite[Theorem 3.1.1(b), p.85]{BK}, the restriction map $H^{0}(X(w_{5}), \mathcal{O}(2k)) \to H^{0}(X(w_{4}), \mathcal{O}(2k))$ is surjective. Further, since $T$ is linearly reductive, the restriction map $H^{0}(X(w_5), \mathcal{O}(2k))^T \to H^{0}(X(w_{4}), \mathcal{O}(2k))^T$ is surjective for all $k \geq 1$. Note that $X_{5}=0$ on $X(w_{4}).$ Consider $X_{i}$'s ($1\le i\le 4$) as elements of $H^0(X(w_{4}),\mathcal{O}(2))^{T}.$ Recall that $(X_1,\ldots,X_4)$ is a basis of $H^0(X(w_{4}),\mathcal{O}(2))^{T}.$
		
		We claim that any relation among $X_{i}$'s $(1\le i\le 4)$ given by a homogeneous polynomials of degree $k$ is identically zero.
		Suppose \begin{equation}\label{eq3.1}
		\sum c_{\underline{m}}X^{\underline{m}}=0
		\end{equation} where $\underline{m}=(m_1,m_2,m_3,m_4)$ are tuples of non-negative integers such that $m_1+m_2+m_3+m_4=k,$ $X^{\underline{m}}=X_1^{m_1}X_2^{m_2}X_3^{m_3}X_4^{m_4}$ and $c_{\underline{m}}$'s are non-zero scalars. Then rewriting \cref{eq3.1} as \begin{equation*}
		\sum_{m_2 \leq m_3}c_{\underline{m}}X_1^{m_1}(X_2X_3)^{m_2}X_3^{m_3-m_2}X_4^{m_4}+\sum_{m_2 >m_3}c_{\underline{m}}X_1^{m_1}X_2^{m_2-m_3}(X_2X_3)^{m_3}X_4^{m_4}=0.
		\end{equation*}
		Recall that by \cref{lemma3.3}, we have $X_2X_3=X_1X_4-Y_1.$ Replacing $X_2X_3$ by $X_1X_4-Y_1$ in the above equation we get 
		\begin{equation*}
		\sum_{m_2 \leq m_3}c_{\underline{m}}X_1^{m_1}(X_1X_4-Y_1)^{m_2}X_3^{m_3-m_2}X_4^{m_4}+\sum_{m_2 >m_3}c_{\underline{m}}X_1^{m_1}X_2^{m_2-m_3}(X_1X_4-Y_1)^{m_3}X_4^{m_4}=0.
		\end{equation*}
	Note that any monomial in $X_1, X_3, X_4, Y_1$ (respectively, in $X_1, X_2,X_4,Y_1$) is standard. Hence, $c_{\underline{m}}=0$ for all $\underline{m}.$ Thus, $X_{1},X_{2},X_{3},X_{4}$ are algebraically independent. Further, by \cref{lemma3.1}, and above surjectivity, the homogeneous coordinate ring of $T\backslash\backslash(X(w_4))^{ss}_T(\mathcal{O}(2))$ is generated by $X_1,X_2,X_3,X_4.$  On the other hand, by \cite[Theorem 3.2.2, p.92]{BK}, $X(w_{4})$ is normal. As $T \backslash \backslash(X(w_{4}))^{ss}_T(\mathcal{O}(2))$ is an open subset of $X(w_{4}),$ it is also normal. Hence, by \cref{rem3.6}, the GIT quotient $T\backslash\backslash(X(w_4))^{ss}_T(\mathcal{O}(2))$ with respect to the descent of $\mathcal{O}(2)$ is projectively normal and isomorphic to $\text{Proj}(\mathbb{ C}[X_{1},X_{2},X_{3}, X_{4}])=\mathbb{P}^{3}$ polarized with $\mathcal{O}(1).$  
		
		By \cite[Theorem 3.1.1(b), p.85]{BK}, it follows that the restriction map $ H^{0}(X(w_{4}), \mathcal{O}(2k)) \to H^{0}(X(w_{i}), \mathcal{O}(2k))$ is surjective for all $1\le i\le 3.$ Further, since $T$ is linearly reductive, the restriction map $ H^{0}(X(w_{4}), \mathcal{O}(2k))^T \to H^{0}(X(w_{i}), \mathcal{O}(2k))^T$ is surjective for all $k \geq 1.$ On the other hand, by \cite[Theorem 3.2.2, p.92]{BK}, $X(w_{i})$ is normal. As $T \backslash \backslash(X(w_{i}))^{ss}_T(\mathcal{O}(2))$ is an open subset of $X(w_{i}),$ it is also normal. Hence, by (iv) the GIT quotient $T\backslash\backslash(X(w_i))^{ss}_T(\mathcal{O}(2))$ with respect to the descent of $\mathcal{O}(2)$ is also projectively normal for all $1\le i\le 3.$
		
		Proof of (iii): Note that $X_{2}$ and $X_{4}$ are identically zero on $X(w_{3}).$ Since the restriction map $H^{0}(X(w_4), \mathcal{O}(2k))^T \to H^{0}(X(w_{3}), \mathcal{O}(2k))^T$ is surjective for all $k \geq 1,$ by (iv) any standard monomial in $H^{0}(X(w_3),\mathcal{O}(2k))^{T}$ is of the form $X_1^{k_1}X_3^{k_2},$ where $k_1+k_2=k.$ Hence, the GIT quotient $T\backslash\backslash(X(w_3))^{ss}_T(\mathcal{O}(2))$ is isomorphic to $\text{Proj}(\mathbb{C}[X_1, X_3])=\mathbb{P}^{1}$ polarized with $\mathcal{O}(1).$
		
		Proof of (ii): Note that $X_{3}$ and $X_{4}$ are identically zero on $X(w_{2}).$ Since the restriction map $H^{0}(X(w_4), \mathcal{O}(2k))^T  \to H^{0}(X(w_{2}), \mathcal{O}(2k))^T$ is surjective for all $k \geq 1,$ by (iv) any standard monomial in $H^{0}(X(w_2),\mathcal{O}(2k))^{T}$ is of the form $X_1^{k_1}X_2^{k_2},$ where $k_1+k_2=k.$ Hence, the GIT quotient $T\backslash\backslash(X(w_2))^{ss}_T(\mathcal{O}(2))$ is isomorphic to $\text{Proj}(\mathbb{C}[X_1, X_2])=\mathbb{P}^{1}$ polarized with $\mathcal{O}(1).$ 
		
		Proof of (i): Note that $X_{2}, X_{3}$ and $X_{4}$ are identically zero on $X(w_{1}).$ Since the restriction map $H^{0}(X(w_4), \mathcal{O}(2k))^T  \to H^{0}(X(w_{1}), \mathcal{O}(2k))^T$ is surjective for all $k \geq 1,$ by (iv) any standard monomial in $H^{0}(X(w_1),\mathcal{O}(2k))^{T}$ is of the form $X_1^k.$ Hence, the GIT quotient $T\backslash\backslash(X(w_1))^{ss}_T(\mathcal{O}(2))$ is isomorphic to $\text{Proj}(\mathbb{C}[X_1]).$ 
		\end{proof}
	
	\begin{remark}
     \item[(i)] The GIT quotient $T\backslash\backslash (G_{1,2})^{ss}_{T}(\mathcal{O}(2))$ with respect to the descent of $\mathcal{O}(2)$ is projectively normal and isomorphic to point.
	
		\item[(ii)] The GIT quotient $T\backslash\backslash (G_{2,4})^{ss}_{T}(\mathcal{O}(2))$ with respect to the descent of $\mathcal{O}(2)$ is projectively normal $($see \cite[Theorem 2.3, p.182]{howard}$)$ and isomorphic to $(\mathbb{P}^{1},\mathcal{O}(1))$ $($see \cite[Proposition 3.5, p.277]{KNS}$)$.
	\end{remark}

Now we prove that the GIT quotient of $X(w_5)$ by $T$ with respect to the descent of $\mathcal{O}(4)$ is projectively normal.

\begin{theorem}\label{thm3.8}
	The homogeneous coordinate ring of $T \backslash \backslash (X(w_5))^{ss}_T(\mathcal{O}(4))$ is generated by elements of degree one.
\end{theorem}
\begin{proof}
Let $f \in H^{0}(X(w_5),\mathcal{O}(4)^{\otimes k})^{T}=H^{0}(X(w_5),\mathcal{O}(4k))^{T}.$ Then by \cref{lemma3.1}, we have $$f=\sum a_{(\underline{m}, \underline{n})}X^{\underline{m}}Y^{\underline{n}},$$ where $\underline{m}=(m_1, m_2, m_3, m_4, m_5),$ $\underline{n}=(n_1, n_2),$ are tuples of non-negative integers such that $m_1+m_2+m_3+m_4+m_5+2n_1+2n_2=2k,$  $X^{\underline{m}}=X_1^{m_1}X_2^{m_2}X_3^{m_3}X_4^{m_4}X_5^{m_5},$   $Y^{\underline{n}}=Y_1^{n_1}Y_2^{n_2},$ and $a_{(\underline{m},\underline{n})}$'s are non-zero scalars.
	
	Now to prove that the homogeneous coordinate ring of $T \backslash \backslash (X(w_5))^{ss}_T(\mathcal{O}(4))$ is generated by $H^{0}(X(w_5),\mathcal{O}(4))^{T}$ as a $\mathbb{C}$-algebra, it is enough to show that for each $f$ as above and each monomial appearing in the expression of $f$ is in the image of $S^
	{k}H^{0}(X(w_5),\mathcal{O}(4))^{T},$ under the natural map $S^
	{k}H^{0}(X(w_5),\mathcal{O}(4))^{T}\to H^0(X(w_5),\mathcal{O}(4k))^{T},$  where  $S^
	{k}H^{0}(X(w_5),\mathcal{O}(4))^{T}$ denotes the $k^{\text{th}}$ symmetric power of the vector space  $H^{0}(X(w_5),\mathcal{O}(4))^{T}.$
	
	Consider the monomial $X^{\underline{m}}Y^{\underline{n}}$ appearing in the expression of $f.$ Since $m_1+m_2+m_3+m_4+m_5$ is even integer, $X^{\underline{m}}$ can be written as $\prod_{(i,j)} X_iX_j,$ where the number of pairs $(i,j)$ is $(k-n_1-n_2)$ and repetition of $X_i$'s are allowed. Thus, $X^{\underline{m}}$ is a product of $(k-n_1-n_2)$ number of monomials in $H^{0}(X(w_5),\mathcal{O}(4))^T.$ Note that $Y_{1},Y_{2}\in H^{0}(X(w_{5}),\mathcal{O}(4))^{T}.$ Therefore, $X^{\underline{m}}Y^{\underline{n}}$ is in the image of  $S^{k}(H^{0}(X(w_5),\mathcal{O}(4))^{T})$ under the natural map.  
\end{proof}

\begin{corollary}\label{cor3.10}
	The GIT quotient $T \backslash \backslash (X(w_5))^{ss}_T(\mathcal{O}(4))$ is projectively normal with respect to the descent of $\mathcal{O}(4).$
\end{corollary}
\begin{proof}
	Follows from \cref{thm3.8}.
\end{proof}

\begin{corollary}
	The GIT quotient $T \backslash \backslash (G_{3,6})^{ss}_T(\mathcal{O}(4))$ is projectively normal with respect to the descent of $\mathcal{O}(4).$
\end{corollary}
\begin{proof}
	Note that for $n=3,$ $w_{5}=(4,5,6).$ So, we have $X(w_{5})=G_{3,6}.$ Therefore, proof immediately follows from \cref{thm3.8}. 	
\end{proof}

In the view of the above results the following question is open:

{\bf Problem:} Is the GIT quotient of $G_{n,2n}$ ($n\ge 4$)  by $T$ with respect to the descent of $\mathcal{O}(4)$ projectively normal?
	
	\section{Appendix}
	Here, we prove the following straightening law on $X(w_5)$ that we used in the proof of \cref{lemma3.3}.
	
	 \tiny{\ytableausetup{boxsize=3.5em} \ytableausetup{centertableaux}\begin{ytableau}
			1 & 3 & 5  & \cdots & 2n-7 & 2n-5 & 2n-2 & 2n-1 \\
			2 & 4 & 6 & \cdots & 2n-6 & 2n-4 & 2n-3 & 2n  
	\end{ytableau}}\normalsize{$ ~=X_1-X_2-X_3+X_4-X_5$}. $(*)$
	\begin{proof}
	Let $\underline{i}=\{1,3,...,2n-7\}$ and $\underline{j}=\{2,4,6,...,2n-6 \}.$
    Let $I=\{1,3,\ldots, 2n-7,2n-5,2n-4\}$ and $J=\{2,4,\ldots,2n-6,2n-3,2n-2,2n-1,2n\}$ be two subsets of $\{1,2,\ldots, 2n\}.$ Then by using (2.1) the following straightening law holds in $X(w_5)$
    
	\noindent
		$p_{\underline{i},2n-5,2n-4,2n-3}p_{\underline{j},2n-2,2n-1,2n}-p_{\underline{i},2n-5,2n-4,2n-2}p_{\underline{j},2n-3,2n-1,2n}+p_{\underline{i},2n-5,2n-4,2n-1}p_{\underline{j},2n-3,2n-2,2n}\\-p_{\underline{i},2n-5,2n-4,2n}p_{\underline{j},2n-3,2n-2,2n-1}=0.$      \hspace{8.7cm}(A1)

	Let $I=\{1,3,\ldots, 2n-7,2n-5,2n-3\}$ and $J=\{2,4,\ldots,2n-6,2n-4,2n-2,2n-1,2n\}.$ 
		Then by using (2.1) the following straightening law holds in  $X(w_5)$
	\noindent	$p_{\underline{i},2n-5,2n-4,2n-3}p_{\underline{j},2n-2,2n-1,2n}+p_{\underline{i},2n-5,2n-3,2n-2}p_{\underline{j},2n-4,2n-1,2n}-p_{\underline{i},2n-5,2n-3,2n-1}p_{\underline{j},2n-4,2n-2,2n}\\+p_{\underline{i},2n-5,2n-3,2n}p_{\underline{j}, 2n-4,2n-2,2n-1}=0.$ \hspace{8.7cm} (A2)
	
	Let $I=\{1,3,\ldots, 2n-7,2n-5,2n-2\}$ and $J=\{2,4,\ldots,2n-6,2n-4,2n-3,2n-1,2n\}.$ Then by using (2.1) the following straightening law holds in $X(w_5)$
	
		\noindent	$p_{\underline{i},2n-5,2n-4,2n-2}p_{\underline{j},2n-3,2n-1,2n}-p_{\underline{i},2n-5,2n-3,2n-2}p_{\underline{j},2n-4,2n-1,2n}-p_{\underline{i},2n-5,2n-2,2n-1}p_{\underline{j},2n-4,2n-3,2n}\\+p_{\underline{i},2n-5,2n-2,2n}p_{\underline{j},2n-4,2n-3,2n-1}=0.$ \hspace{8.7cm} (A3)
	
	Let $I=\{1,3,\ldots, 2n-7,2n-5,2n-1\}$ and $J=\{2,4,\ldots,2n-6,2n-4,2n-3,2n-2,2n\}.$ Then by using (2.1) the following straightening law holds in $X(w_5)$
	
	\noindent	$p_{\underline{i},2n-5,2n-4,2n-1}p_{\underline{j},2n-3,2n-2,2n}-p_{\underline{i},2n-5,2n-3,2n-1}p_{\underline{j},2n-4,2n-2,2n}+p_{\underline{i},2n-5,2n-2,2n-1}p_{\underline{j},2n-4,2n-3,2n}\\+p_{\underline{i},2n-5,2n-1,2n}p_{\underline{j},2n-4,2n-3,2n-2}=0.$ \hspace{8.7cm} (A4)
	
	Let $I=\{1,3,\ldots, 2n-7,2n-5,2n\}$ and $J=\{2,4,\ldots,2n-6,2n-4,2n-3,2n-2,2n-1\}.$ Then by using (2.1) the following straightening law holds in $X(w_5)$
	
	\noindent	$p_{\underline{i},2n-5,2n-4,2n}p_{\underline{j},2n-3,2n-2,2n-1}-p_{\underline{i},2n-5,2n-3,2n}p_{\underline{j},2n-4,2n-2,2n-1}+p_{\underline{i},2n-5,2n-2,2n}p_{\underline{j},2n-4,2n-3,2n-1}\\-p_{\underline{i},2n-5,2n-1,2n}p_{\underline{j},2n-4,2n-3,2n-2}=0.$ \hspace{8.7cm}(A5)
	
	Thus by using (A3), we have
	
	\noindent	$p_{\underline{i},2n-5,2n-2,2n-1}p_{\underline{j},2n-4,2n-3,2n}=p_{\underline{i},2n-5,2n-4,2n-2}p_{\underline{j},2n-3,2n-1,2n}~-~p_{\underline{i},2n-5,2n-3,2n-2}p_{\underline{j},2n-4,2n-1,2n}$\\$~+~p_{\underline{i},2n-5,2n-2,2n}p_{\underline{j},2n-4,2n-3,2n-1}.$
	
	By using (A5) we have
	
	\noindent	$p_{\underline{i},2n-5,2n-2,2n-1}p_{\underline{j},2n-4,2n-3,2n}=p_{\underline{i},2n-5,2n-4,2n-2}p_{\underline{j},2n-3,2n-1,2n}-p_{\underline{i},2n-5,2n-3,2n-2}p_{\underline{j},2n-4,2n-1,2n}
		-p_{\underline{i},2n-5,2n-4,2n}p_{\underline{j},2n-3,2n-2,2n-1}+p_{\underline{i},2n-5,2n-3,2n}p_{\underline{j},2n-4,2n-2,2n-1}+p_{\underline{i},2n-5,2n-1,2n}p_{\underline{j},2n-4,2n-3,2n-2}.$
	
Further, by using (A1), (A2) and (A4),  (*) follows. 
\end{proof}	
{\bf Acknowledgements.} We are grateful to the referee for careful reading, numerous valuable comments and the reference of the article  C. Teleman,  The quantization conjecture revisited, which helps to improve the exposition of this article.

\end{document}